\documentclass{amsart}
\usepackage{amssymb}
\usepackage{amsmath}
\usepackage{amsfonts}

\newtheorem{theorem}{Theorem}
\theoremstyle{plain}

\numberwithin{equation}{section}
\input{tcilatex}

\begin{document}
\title[Ricci solitons on LBCV spaces]{Ricci solitons of three-dimensional
Lorentzian Bianchi-Cartan-Vranceanu spaces}
\author{Murat ALTUNBA\c{S}}
\address{Erzincan Binali Y\i ld\i r\i m University, Faculty of Science and
Art, Department of Mathematics, 24030, Erzincan-Turkey.}
\email{maltunbas@erzincan.edu.tr}

\begin{abstract}
Explicit formulae for homogenous Ricci solitons on three-dimensional
Lorentzian Bianchi-Cartan-Vranceanu spaces are obtained.

\textbf{AMS Classification [2020]: }53C30, 53C50.

\end{abstract}

\maketitle

\section{Introduction}

A Ricci soliton metric on a manifold $M$ is defined by the condition%
\begin{equation}
L_{X}g+\rho =\gamma g,  \label{31}
\end{equation}%
where $X$ is a smooth vector field on $M$, $L_{X}g$ is Lie derivative in the
direction of $X$ and $\gamma $ is a real constant. A Ricci soliton is called 
\textit{shrinking} if $\gamma >0$, \textit{steady} if $\gamma =0$ and 
\textit{expanding} if $\gamma <0.$ Ricci soliton metrics are a
generalization of Einstein metrics.

Ricci solitons have been extensively studied in many works from many points
of view, so we refer only \cite{Cao} and \cite{Chow} for more information
about geometry of Ricci solitons.

Many researchers have been particularly interested in Ricci solitons on
three-dimensional homogenous spaces, such as the Lie group $SL(2,%
\mathbb{R}
),\ $Heisenberg group $Nil_{3},$ Berger spheres $S_{Berger}^{3},$ $%
S^{2}\times 
\mathbb{R}
$, $H^{2}\times 
\mathbb{R}
$ and the Lorentzian-Heisenberg group (see \cite{Baird},\cite{BatatOnda},%
\cite{Jablonski},\cite{Onda},\cite{Vazquez}).

Bianchi-Cartan-Vranceanu spaces are three-dimensional homogenous spaces with
four dimensional isometry group. Ricci solitons on Bianchi-Cartan-Vranceanu
spaces were studied by Batat \textit{et al. }in \cite{Batat}.

Lorentzian Bianchi-Cartan-Vranceanu spaces (briefly LBCV-spaces) are
considered by several authors in very recent papers, especially when
investigating some special curves such as slant, Legendre and biharmonic
etc. on it (see \cite{Lee1}, \cite{Lee2}, \cite{Yildirim}).

As we mentioned above, although the subject of Ricci solitons is
well-studied on homogenous manifolds, we give a classification of Ricci
solitons by obtaining explicit formulae on LBCV-spaces in this paper. In
fact, we will prove the following theorem:

\begin{theorem}
\label{Theorem} Let LBCV-spaces with the metric in (\ref{32}) are given.
Then the following statements are true:

(i) LBCV-spaces do not admit homogenous Ricci solitons when $\lambda \neq 0$
and $\mu >0.$

(ii) LBCV-spaces admit shrinking homogenous Ricci solitons when $\lambda
\neq 0$ and $\mu =0.$

(iii) LBCV-spaces admit expanding homogenous Ricci solitons when $\lambda
\neq 0$ and $\mu <0.$

(iv) LBCV-spaces admit shrinking homogenous Ricci solitons when $\lambda =0$
and $\mu >0.$

(v) LBCV-spaces admit expanding homogenous Ricci solitons when $\lambda =0$
and $\mu <0.$
\end{theorem}

\section{\protect\bigskip Lorentzian Bianchi-Cartan-Vranceanu spaces
(LBCV-spaces)}

In this section, we will recall some fundamental properties of LBCV-spaces
(see \cite{Yildirim}, \cite{Lee1}).

Let $\lambda ,\mu \in 
\mathbb{R}
.$ An open subset of $%
\mathbb{R}
^{3}$ is given by

\begin{equation*}
D=\{(x,y,z)\in 
\mathbb{R}
^{3}:1+\mu (x^{2}+y^{2})>0\}.
\end{equation*}%
The Lorentzian metric is equipped as following:%
\begin{equation}
g_{\lambda ,\mu }=\frac{dx^{2}+dy^{2}}{(1+\mu (x^{2}+y^{2}))^{2}}-\left( dz+%
\frac{\lambda }{2}\frac{ydx-xdy}{1+\mu (x^{2}+y^{2})}\right) ^{2}.
\label{32}
\end{equation}%
The pair $(D,g_{\lambda ,\mu })$ is called Lorentzian
Bianchi-Cartan-Vranceanu spaces and it is denoted by $M_{\lambda ,\mu }.$

An orthonormal frame field is given by 
\begin{equation}
E_{1}=\delta \frac{\partial }{\partial x}-\frac{\lambda y}{2}\frac{\partial 
}{\partial z},\ E_{2}=\delta \frac{\partial }{\partial y}+\frac{\lambda x}{2}%
\frac{\partial }{\partial z},\ E_{3}=\frac{\partial }{\partial z},
\label{33}
\end{equation}%
where we write $\delta =1+\mu (x^{2}+y^{2}).$

Therefore, the Lie brackets are obtained as%
\begin{equation*}
\lbrack E_{1},E_{2}]=-2\mu yE_{1}+2\mu xE_{2}+\lambda E_{3},\
[E_{1},E_{3}]=[E_{2},E_{3}]=0.
\end{equation*}%
Let $\nabla $ and $R$ denote the Levi-Civita connection and the curvature
tensor of $M_{\lambda ,\mu },$ respectively. We have 
\begin{eqnarray*}
\nabla _{E_{1}}E_{1} &=&2\mu yE_{2},\ \nabla _{E_{1}}E_{2}=-2\mu yE_{1}+%
\frac{\lambda }{2}E_{3},\ \nabla _{E_{1}}E_{3}=\frac{\lambda }{2}E_{2}, \\
\nabla _{E_{2}}E_{1} &=&-2\mu xE_{2}-\frac{\lambda }{2}E_{3},\ \nabla
_{E_{2}}E_{2}=2\mu xE_{1},\ \nabla _{E_{2}}E_{3}=-\frac{\lambda }{2}E_{1}, \\
\nabla _{E_{3}}E_{1} &=&\frac{\lambda }{2}E_{2},\ \nabla _{E_{3}}E_{2}=-%
\frac{\lambda }{2}E_{1},\ \nabla _{E_{3}}E_{3}=0.
\end{eqnarray*}%
The components of the curvature tensor $R_{ijk}^{l}$ are given by \cite%
{Yildirim2} 
\begin{eqnarray*}
R_{121}^{1} &=&0,\ R_{313}^{1}=\frac{\lambda ^{2}}{4},\ R_{323}^{1}=0,\
R_{221}^{1}=-4\mu -\frac{3}{4}\lambda ^{2},\ R_{331}^{1}=-\frac{\lambda ^{2}%
}{4}, \\
R_{112}^{1} &=&0,\ R_{223}^{1}=0,\ R_{212}^{1}=4\mu +\frac{3}{4}\lambda
^{2},\ R_{332}^{1}=0,\ R_{113}^{1}=0,
\end{eqnarray*}%
\begin{eqnarray*}
R_{121}^{2} &=&4\mu +\frac{3}{4}\lambda ^{2},\ R_{313}^{2}=0,\ R_{323}^{2}=%
\frac{\lambda ^{2}}{4},\ R_{221}^{2}=0,\ R_{331}^{2}=0, \\
R_{112}^{2} &=&-4\mu -\frac{3}{4}\lambda ^{2},\ R_{223}^{2}=0,\
R_{212}^{2}=0,\ R_{332}^{2}=-\frac{\lambda ^{2}}{4},\ R_{113}^{2}=0,
\end{eqnarray*}%
\begin{eqnarray*}
R_{121}^{3} &=&0,\ R_{313}^{3}=0,\ R_{323}^{3}=0,\ R_{221}^{3}=0,\
R_{331}^{3}=0, \\
R_{112}^{3} &=&0,\ R_{223}^{3}=-\frac{\lambda ^{2}}{4},\ R_{212}^{3}=0,\
R_{332}^{3}=0,\ R_{113}^{3}=-\frac{\lambda ^{2}}{4}.
\end{eqnarray*}%
Therefore, for the Ricci tensor $\rho (X,Y)=tr\{Z\rightarrow R(X,Z)Y\}$ with
respect to orthonormal basis (\ref{33}), we obtain%
\begin{equation}
\rho _{11}=\rho _{22}=4\mu +\lambda ^{2},\ \rho _{33}=\frac{\lambda ^{2}}{2},
\label{34}
\end{equation}%
where we set $\rho _{ij}=\rho (E_{i},E_{j}).$

\section{\protect\bigskip Ricci solitons on Lorentzian
Bianchi-Cartan-Vranceanu spaces}

In this section we deal with the Ricci solitons on LBCV-space $M_{\lambda
,\mu }=(D,g_{\lambda ,\mu }).$ Let $X=X_{1}E_{1}+X_{2}E_{2}+X_{3}E_{3}$ be
an arbitrary vector field on $M_{\lambda ,\mu },$ where $X_{1},X_{2},X_{3}$
are smooth functions of the variables $x,y,z.$ Then, the Lie derivative of
the metric (\ref{32}) satisfies the following relations:

\begin{eqnarray}
L_{X}g_{\lambda ,\mu }(E_{1},E_{1}) &=&2(E_{1}(X_{1})-2\mu yX_{2}),
\label{35} \\
L_{X}g_{\lambda ,\mu }(E_{1},E_{2}) &=&2\mu xX_{2}+2\mu
yX_{1}+E_{1}(X_{2})+E_{2}(X_{1}),  \notag \\
L_{X}g_{\lambda ,\mu }(E_{1},E_{3}) &=&E_{3}(X_{1})-E_{1}(X_{3})-\lambda
X_{2},  \notag \\
L_{X}g_{\lambda ,\mu }(E_{2},E_{2}) &=&2(E_{2}(X_{2})-2\mu xX_{1}),  \notag
\\
L_{X}g_{\lambda ,\mu }(E_{2},E_{3}) &=&\lambda
X_{1}-E_{2}(X_{3})+E_{3}(X_{2}),  \notag \\
L_{X}g_{\lambda ,\mu }(E_{3},E_{3}) &=&-2E_{3}(X_{3}).  \notag
\end{eqnarray}%
Therefore, if we use (\ref{32}), (\ref{34}) and (\ref{35}) in (\ref{31}) and
have in mind (\ref{33}), with a standard calculation, we see that a LBCV
space is Ricci soliton if and only if the following system is satisfied: 
\begin{equation}
\begin{array}{c}
2\mu yX_{2}-\delta \partial _{x}X_{1}+\frac{\lambda }{2}y\partial _{z}X_{1}=%
\frac{\rho _{11}-\gamma }{2}, \\ 
2\mu xX_{2}+2\mu yX_{1}+\delta \partial _{x}X_{2}-\frac{\lambda }{2}%
y\partial _{z}X_{2}+\delta \partial _{y}X_{1}+\frac{\lambda }{2}x\partial
_{z}X_{1}=0, \\ 
-\lambda X_{2}-\delta \partial _{x}X_{3}+\frac{\lambda }{2}y\partial
_{z}X_{3}+\partial _{z}X_{1}=0, \\ 
2\mu xX_{1}-\delta \partial _{y}X_{2}-\frac{\lambda }{2}x\partial _{z}X_{2}=%
\frac{\rho _{11}-\gamma }{2}, \\ 
\lambda X_{1}-\delta \partial _{y}X_{3}-\frac{\lambda }{2}x\partial
_{z}X_{3}+\partial _{z}X_{2}=0, \\ 
\partial _{z}X_{3}=\frac{\gamma +\rho _{33}}{2},%
\end{array}
\label{36}
\end{equation}%
where we set $\partial _{x}=\frac{\partial }{\partial x},\partial _{y}=\frac{%
\partial }{\partial y},\partial _{z}=\frac{\partial }{\partial z}.$

Equation (\ref{36})$_{6}$ implies that 
\begin{equation}
X_{3}=(\frac{\gamma +\rho _{33}}{2})z+A(x,y),\ A\in C^{\infty }(M).
\label{37}
\end{equation}%
for an arbitrary smooth function $A=A(x,y).$

\subsection{Case 1 $\protect\lambda \neq 0$}

From (\ref{36})$_{5}$ and using (\ref{37}), we get

\begin{equation}
X_{1}=\frac{1}{\lambda }\left( \delta \partial _{y}A-\partial
_{z}X_{2}+\lambda (\frac{\gamma +\rho _{33}}{4})x\right) .  \label{38}
\end{equation}%
Substituting (\ref{37}) and (\ref{38}) in (\ref{36})$_{3},$ we occur 
\begin{equation*}
\lambda ^{2}X_{2}+\partial _{z}^{2}X_{2}=\lambda \left( \lambda (\frac{%
\gamma +\rho _{33}}{4})y-\delta \partial _{x}A\right) .
\end{equation*}%
Solution of the above equation gives us%
\begin{equation}
X_{2}=-\frac{\delta }{\lambda }\partial _{x}A+(\frac{\gamma +\rho _{33}}{4}%
)y+C_{1}(x,y)\cos (\lambda z)+C_{2}(x,y)\sin (\lambda z),  \label{39}
\end{equation}%
where $C_{1}$ and $C_{2}$ are arbitrary smooth functions of the variables $x$
and $y.$

It follows that 
\begin{equation}
X_{1}=\frac{\delta }{\lambda }\partial _{y}A+(\frac{\gamma +\rho _{33}}{4}%
)x+C_{1}(x,y)\sin (\lambda z)-C_{2}(x,y)\cos (\lambda z).  \label{40}
\end{equation}%
By substituting (\ref{39}) and (\ref{40}) in (\ref{36})$_{1},\ $we see that 
\begin{equation}
\begin{array}{c}
\partial _{x}C_{1}=\left( 2\mu +\frac{\lambda ^{2}}{2}\right) \frac{yC_{2}}{%
\delta }, \\ 
\partial _{x}C_{2}=-\left( 2\mu +\frac{\lambda ^{2}}{2}\right) \frac{yC_{1}}{%
\delta }, \\ 
(1+\mu (x^{2}-y^{2}))(\frac{\gamma +\rho _{33}}{4})+\frac{\delta }{\lambda }%
\left( 2\mu (x\partial _{y}A+y\partial _{x}A)+\delta \partial _{x}\partial
_{y}A\right) =\frac{\gamma -\rho _{11}}{2}{.}%
\end{array}
\label{41}
\end{equation}%
Again, by substituting (\ref{39}) and (\ref{40}) in (\ref{36})$_{4},\ $we
obtain%
\begin{equation}
\begin{array}{c}
\partial _{y}C_{1}=-\left( 2\mu +\frac{\lambda ^{2}}{2}\right) \frac{xC_{2}}{%
\delta }, \\ 
\partial _{y}C_{2}=\left( 2\mu +\frac{\lambda ^{2}}{2}\right) \frac{xC_{1}}{%
\delta }, \\ 
(1-\mu (x^{2}-y^{2}))(\frac{\gamma +\rho _{33}}{4})-\frac{\delta }{\lambda }%
\left( 2\mu (x\partial _{y}A+y\partial _{x}A)+\delta \partial _{x}\partial
_{y}A\right) =\frac{\gamma -\rho _{11}}{2}{.}%
\end{array}
\label{42}
\end{equation}%
The last equations in (\ref{41}) and (\ref{42}) show that 
\begin{eqnarray*}
\gamma  &=&2\rho _{11}+\rho _{33} \\
\gamma  &=&8\mu +\frac{3\lambda ^{2}}{2}.
\end{eqnarray*}%
Therefore, (\ref{41}) and (\ref{42}) turn to be 
\begin{equation}
\lambda \mu \left( 2\mu +\frac{\lambda ^{2}}{2}\right) (x^{2}-y^{2})+\delta
\left( 2\mu (x\partial _{y}A+y\partial _{x}A)+\delta \partial _{x}\partial
_{y}A\right) =0.  \label{43}
\end{equation}%
Taking derivative with respect to $y$ in the first equation of (\ref{41})
and with respect to $x$ in the first equation of (\ref{42}), and having in
mind $\partial _{x}C_{2}$ and $\partial _{y}C_{2},$ we see that $C_{2}=0$
(when $\lambda ^{2}\neq -4\mu )$ or $C_{2}\in 
\mathbb{R}
$ (when $\lambda ^{2}=-4\mu ).$ Similarly, $C_{1}$ is zero or constant.

Let the inequality $\lambda ^{2}\neq -4\mu $ holds. Equation (\ref{36})$_{2}$
leads to 
\begin{equation}
2\lambda \mu \left( 4\mu +\lambda ^{2}\right) xy+\delta \lbrack \left( 4\mu
(y\partial _{y}A-x\partial _{x}A)+\delta (\partial _{y}^{2}A-\partial
_{x}^{2}A\right) ]=0.  \label{44}
\end{equation}%
So, the vector field $X=X_{1}E_{1}+X_{2}E_{2}+X_{3}E_{3}$ fulfils (\ref{36})
if and only if 
\begin{equation*}
\begin{array}{c}
X_{1}=\frac{\delta }{\lambda }\partial _{y}A+\left( \frac{4\mu +\lambda ^{2}%
}{2}\right) x, \\ 
X_{2}=-\frac{\delta }{\lambda }\partial _{x}A+\left( \frac{4\mu +\lambda ^{2}%
}{2}\right) y, \\ 
X_{3}=(4\mu +\lambda ^{2})z+A{.}%
\end{array}%
\end{equation*}%
Here, the function $A$ satisfies (\ref{43}) and (\ref{44}).

Now, suppose that $\lambda ^{2}=-4\mu .$ In this case, Equations (\ref{43})
and (\ref{44}) remain valid, but the vector field $X$ reduces to 
\begin{equation}
\begin{array}{c}
X_{1}=\frac{\delta }{\lambda }\partial _{y}A+C_{1}\sin (\lambda z)-C_{2}\cos
(\lambda z), \\ 
X_{2}=-\frac{\delta }{\lambda }\partial _{x}A+C_{1}\cos (\lambda
z)+C_{2}\sin (\lambda z), \\ 
X_{3}=A{,}%
\end{array}
\label{45}
\end{equation}%
$C_{1},C_{2}\in 
\mathbb{R}
$ and $\gamma =2\mu .$

(a) If $\mu =0$, Equations (\ref{43}) and (\ref{44}) turn in to be 
\begin{equation*}
\partial _{x}\partial _{y}A=0\ \text{and }\partial _{y}^{2}A=\partial
_{x}^{2}A.
\end{equation*}%
So, we have%
\begin{equation*}
A=a_{1}(x^{2}+y^{2})+a_{2}x+a_{3}y+a_{4},\ a_{1},...,a_{4}\in 
\mathbb{R}
.
\end{equation*}%
As a result, when $\mu =0,$ the vector field $%
X=X_{1}E_{1}+X_{2}E_{2}+X_{3}E_{3}$ satisfy the soliton equation (\ref{31})
if and only if 
\begin{equation*}
\begin{array}{c}
X_{1}=\frac{1}{\lambda }(2a_{1}y+a_{3})-\frac{\lambda ^{2}}{4}x, \\ 
X_{2}=-\frac{1}{\lambda }(2a_{1}x+a_{2})-\frac{\lambda ^{2}}{4}y, \\ 
X_{3}=-\frac{\lambda ^{2}}{2}z+a_{1}(x^{2}+y^{2})+a_{2}x+a_{3}y+a_{4}{,}%
\end{array}%
\end{equation*}%
where $a_{1},...,a_{4}\in 
\mathbb{R}
$ and $\gamma =\frac{3\lambda ^{2}}{2}>0.$ Thus, we proved Theorem \ref%
{Theorem} (ii).

(b) Now, suppose that $\mu \neq 0$. Set $f=\delta A$ and $\Delta =\lambda
\mu \left( 2\mu +\frac{\lambda ^{2}}{2}\right) .$ Then Equations (\ref{43})
and (\ref{44}) imply 
\begin{equation}
\partial _{x}\partial _{y}f=\frac{\Delta (y^{2}-x^{2})}{1+\mu (x^{2}+y^{2})},
\label{46}
\end{equation}%
\begin{equation}
\partial _{x}^{2}f-\partial _{y}^{2}f=\frac{4\Delta xy}{1+\mu (x^{2}+y^{2})}.
\label{47}
\end{equation}%
If we integrate (\ref{46}) with respect to $y$, we get 
\begin{equation}
\partial _{x}f=\Delta \left[ \frac{y}{\mu }-\frac{(1+2\mu x^{2})}{\left\vert
\mu \right\vert ^{3/2}\sqrt{1+\mu x^{2}}}\arctan \left( \frac{\sqrt{%
\left\vert \mu \right\vert y}}{\sqrt{1+\mu x^{2}}}\right) \right] +\alpha
(x),  \label{48}
\end{equation}%
and if we integrate (\ref{46}) with respect to $x$, we obtain 
\begin{equation}
\partial _{y}f=\Delta \left[ -\frac{x}{\mu }+\frac{(1+2\mu y^{2})}{%
\left\vert \mu \right\vert ^{3/2}\sqrt{1+\mu y^{2}}}\arctan \left( \frac{%
\sqrt{\left\vert \mu \right\vert x}}{\sqrt{1+\mu y^{2}}}\right) \right]
+\beta (y),  \label{49}
\end{equation}%
where $\alpha $ and $\beta $ are smooth functions. Remark that if $\mu <0,$
we have $arctanh$ instead of $arctan.$ Differentiating (\ref{48}) by $x$ and
(\ref{49}) by $y$, replacing into (\ref{47}), we deduce that there is a
solution if and only if $\Delta =0,$ that is, if $\mu =-\frac{\lambda ^{2}}{4%
}<0.\ $This shows that when $\mu >0$ the solution does not exist which
proves the statement Theorem \ref{Theorem} (i). Moreover, we occur that 
\begin{equation*}
f=a_{1}(x^{2}+y^{2})+a_{2}x+a_{3}y+a_{4},
\end{equation*}%
\begin{equation*}
\text{and }A(x,y)=\frac{a_{1}(x^{2}+y^{2})+a_{2}x+a_{3}y+a_{4}}{1+\mu
(x^{2}+y^{2})}.
\end{equation*}%
Thus, if $\mu >0$ Equation (\ref{31}) has no solution and if $\mu <0$ it is
satisfied only for $\mu =-\frac{\lambda ^{2}}{4}.$ Then, from (\ref{45}), we
obtain the corresponding solutions as follows: 
\begin{equation}
\begin{array}{c}
X_{1}=\frac{-2a_{2}\mu xy+a_{3}\left( \mu (x^{2}-y^{2})+1\right) -2a_{4}\mu
y+2a_{1}y}{\lambda (1+\mu (x^{2}+y^{2}))} \\ 
+a_{5}\sin (\lambda z)-a_{6}\cos (\lambda z), \\ 
X_{2}=\frac{2\mu x(a_{3}y+a_{4})+a_{2}\left( \mu (x^{2}-y^{2})-1\right)
-2a_{1}x}{\lambda (1+\mu (x^{2}+y^{2}))} \\ 
+a_{5}\cos (\lambda z)+a_{6}\sin (\lambda z), \\ 
X_{3}=\frac{a_{1}(x^{2}+y^{2})+a_{2}x+a_{3}y+a_{4}}{1+\mu (x^{2}+y^{2})}{,}%
\end{array}
\label{50}
\end{equation}%
with $a_{1},...,a_{6}\in 
\mathbb{R}
$ and $\gamma =-\frac{\lambda ^{2}}{2}<0.$ This completes the proof of
Theorem \ref{Theorem} (iii). Remark that in this case associated the
solitons are Killing vector fields also.

\subsection{Case 2 $\protect\lambda =0,\ \protect\mu \neq 0$}

In this case the system (\ref{36}) reduces to

\begin{equation}
\begin{array}{c}
2\mu yX_{2}-\delta \partial _{x}X_{1}=\frac{4\mu -\gamma }{2}, \\ 
2\mu xX_{2}+2\mu yX_{1}+\delta \partial _{x}X_{2}+\delta \partial
_{y}X_{1}=0, \\ 
-\delta \partial _{x}X_{3}+\partial _{z}X_{1}=0, \\ 
2\mu xX_{1}-\delta \partial _{y}X_{2}=\frac{4\mu -\gamma }{2}, \\ 
-\delta \partial _{y}X_{3}+\partial _{z}X_{2}=0, \\ 
\partial _{z}X_{3}=\frac{\gamma }{2}.%
\end{array}
\label{51}
\end{equation}%
From the equations (\ref{51})$_{3}$, (\ref{51})$_{5}$ and (\ref{51})$_{6}$,
we obtain%
\begin{equation}
\begin{array}{c}
X_{1}=\delta (\partial _{x}A)z+F(x,y), \\ 
X_{2}=\delta (\partial _{y}A)z+E(x,y), \\ 
X_{3}=\frac{\gamma }{2}z+A(x,y){,}%
\end{array}
\label{52}
\end{equation}%
where $A,E$ and $F$ are smooth functions of $x$ and $y.$ Putting these
expressions of $X_{1}$ and $X_{2}$ in (\ref{51})$_{1}$ gives us 
\begin{equation*}
-\delta \lbrack 2\mu ((x\partial _{x}A-y\partial _{y}A)+\delta \partial
_{x}^{2}A]z+2\mu yE-\delta \partial xF=\frac{4\mu -\gamma }{2}.
\end{equation*}%
Since this equation holds for all $z,$ we have%
\begin{equation}
2\mu ((x\partial _{x}A-y\partial _{y}A)+\delta \partial _{x}^{2}A=0,\ 2\mu
yE-\delta \partial xF=\frac{4\mu -\gamma }{2}.  \label{53}
\end{equation}%
Again, substituting the expressions of $X_{1}$ and $X_{2}$ in (\ref{52})
into (\ref{51})$_{4}$ and (\ref{51})$_{2}$ we obtain, respectively 
\begin{equation}
2\mu ((y\partial _{y}A-x\partial _{x}A)+\delta \partial _{y}^{2}A=0,\ 2\mu
xF-\delta \partial yE=\frac{4\mu -\gamma }{2},  \label{54}
\end{equation}%
and%
\begin{equation}
2\mu (x\partial _{y}A+y\partial _{x}A)+\delta \partial _{x}\partial
_{y}A=0,\ 2\mu (xE+yF)+\delta (\partial _{x}E+\partial _{y}F)=0.  \label{55}
\end{equation}%
Combining the first equations in (\ref{53}) and (\ref{54}), we get 
\begin{equation}
\partial _{x}^{2}A+\partial _{y}^{2}A=0.  \label{56}
\end{equation}%
If we derive the first equation in (\ref{53}) with respect to $x$ and the
first equation with respect to $y$ in (\ref{54}), and have in mind (\ref{56}%
), we occur 
\begin{equation}
2\partial _{x}A+x\partial _{x}^{2}A+y\partial _{x}\partial _{y}A=0.
\label{57}
\end{equation}%
Now, if we derive the first equation in (\ref{53}) with respect to $y$ and
the first equation with respect to $x$ in (\ref{54}), and by virtue of (\ref%
{56}), we deduce 
\begin{equation}
2\partial _{y}A-y\partial _{x}^{2}A+x\partial _{x}\partial _{y}A=0.
\label{58}
\end{equation}%
Therefore, from (\ref{57}) and (\ref{58}), after using the first equation in
(\ref{53}), we obtain that $\partial _{x}^{2}A=\partial _{y}^{2}A=0.$ So,
the first equations in (\ref{53}) and (\ref{54}) become $x\partial
_{x}A-y\partial _{y}A=0,$ which together with the first equation of \ (\ref%
{55}) shows that $A$ is a constant function.

Similarly, by considering the second equations of (\ref{53}), (\ref{54}) and
(\ref{55}), we have 
\begin{equation*}
\partial _{y}(\delta E)-\partial _{x}(\delta F)=0,\ \partial _{x}(\delta
E)+\partial _{y}(\delta F)=0.
\end{equation*}%
The solution of this system is $\delta E=c_{1}$, $\delta F=c_{2},$ where $%
c_{1},c_{2}\in 
\mathbb{R}
.$ Putting this in (\ref{53}), we obtain $E=F=0$ with $\gamma =4\mu .$ So,
by setting $A=a\in 
\mathbb{R}
,$ the system (\ref{52}) turns in to be \ 
\begin{equation*}
X_{1}=X_{2}=0,\ X_{3}=2\mu z+a.
\end{equation*}%
This completes the proof of Theorem \ref{Theorem} 1.

\subsection{Case 3 $\protect\lambda =\ \protect\mu =0.$}

In this final case we deal with the Minkowski three-space. If $\lambda =\
\mu =0,$ the system (\ref{52}) becomes

\begin{equation*}
\begin{array}{c}
\delta \partial _{x}X_{1}=\frac{\gamma }{2}, \\ 
\partial _{x}X_{2}+\partial _{y}X_{1}=0, \\ 
-\partial _{x}X_{3}+\partial _{z}X_{1}=0, \\ 
\partial _{y}X_{2}=\frac{\gamma }{2}, \\ 
-\partial _{y}X_{3}+\partial _{z}X_{2}=0, \\ 
\partial _{z}X_{3}=\frac{\gamma }{2}.%
\end{array}%
\end{equation*}%
By direct computation, we see that, for $X=X_{1}E_{1}+X_{2}E_{2}+X_{3}E_{3},%
\ $the corresponding soliton has the following form:%
\begin{equation*}
\begin{array}{c}
X_{1}=\frac{\gamma }{2}x-a_{1}y+a_{2}z+a_{3}, \\ 
X_{2}=a_{1}x+\frac{\gamma }{2}y+a_{4}z+a_{5}, \\ 
X_{3}=a_{2}x+a_{4}y+\frac{\gamma }{2}z+a_{6}{,}%
\end{array}%
\end{equation*}%
for every $\gamma \in 
\mathbb{R}
$ with $a_{1},...,a_{6}\in 
\mathbb{R}
.$

\section{Conclusion}

In this work, we gave a classification for Ricci solitons on Lorentzian
Bianchi-Cartan-Vranceanu spaces. We showed that there exist significant
differences from the Riemannian case, which is studied in the reference \cite%
{Batat}, when $\lambda \neq 0.$


\begin{thebibliography}{10}
\bibitem[1]{Baird} P. Baird and L. Danielo, \emph{Three-dimensional Ricci
solitons which project to surfaces, }J. Reine Angew. Math. \textbf{608}
(2007), 65---91.

\bibitem[2]{Batat} W. Batat, T. Sukilovic and S. Vukmirovic, \emph{Ricci
solitons of three-dimensional Bianchi-Cartan-Vranceanu spaces}, J. Geom., 
\textbf{111:1} (2020), 1-10.

\bibitem[3]{BatatOnda} W. Batat and K. Onda,\emph{\ Algebraic Ricci solitons
of three-dimensional Lorentzian Lie groups.} Journal of Geometry and
Physics, \textbf{114} (2017) 138--152.

\bibitem[4]{Cao} H. D. Cao, \emph{Geometry of Ricci solitons}, Chinese Ann.
Math. Ser. B \textbf{27B }(2006), 121-142.

\bibitem[5]{Chow} B. Chow, D. Knopf, \emph{The Ricci Flow: An Introduction,
Mathematical Surveys and Monographs}, \textbf{110}. American Mathematical
Society, Providence, 2004.

\bibitem[6]{Jablonski} M. Jablonski, \emph{Homogenous Ricci solitons}, J.
Reine Angew. Math.\textbf{\ 699} (2015), 159-182.

\bibitem[7]{Lee1} J.E. Lee, \emph{Slant curves in contact Lorentzian
manifolds with CR structures}. Mathematics. \textbf{8 }(1) (2020), 46.

\bibitem[8]{Lee2} J.E. Lee, \emph{Biharmonic curves in 3-dimensional
Lorentzian-Sasakian space forms}, Comm. Korean Math. Soc., \textbf{35 }(3)%
\textbf{\ }(2020), 967-977.

\bibitem[9]{Onda} K. Onda. \emph{Lorentz Ricci Solitons on 3-dimensional Lie
groups.} Geom Dedicata \textbf{147} (2010), 313--322.

\bibitem[10]{Vazquez} M. B. Vazquez, G. Calvaruso, E. Garc\'{\i}a-R\'{\i}o
and S. Gavino-Fern\'{a}ndez, \emph{Three-dimensional Lorentzian homogeneous
Ricci solitons}, Isr. J. Math. \textbf{188} (2012), 385--403.

\bibitem[11]{Yildirim} A. Yildirim, \emph{Slant curve in Lorentzian BCV
spaces, }J.\emph{\ }Geometry and Symmetry in Physics, \textbf{56} (2020),
67-85.

\bibitem[12]{Yildirim2} A.Yildirim, On Lorentzian BCV spaces, Int. J. Math.
Archive, \textbf{3 }(4) (2012), 1365-1371.
\end{thebibliography}
\end{document}